\documentclass{tac}
\usepackage{amsmath}
\usepackage{amssymb}
\usepackage{stmaryrd}
\usepackage[all,cmtip,2cell]{xy}
\usepackage{bibentry}

\author{Leonardo Larizza}

\thanks{The author gratefully acknowledges financial support from the Centre for Mathematics
of the University of Coimbra (UIDB/00324/2020, funded by the Portuguese Government through
FCT/MCTES).
The author also acknowledges the support of a PhD grant from FCT/MCTES
(PD/BD/142958/2018).
This work was developed in the context of his PhD thesis work under the supervision of Maria Manuel Clementino.}

\address{University of Coimbra, CMUC, Department of Mathematics, 3001-501 Coimbra, Portugal}

\title[On factorisations for Ord-enriched categories and partial maps]{On factorisation systems for Ord-enriched categories and categories of partial maps}

\copyrightyear{2021}



\eaddress{leonardo.larizza91@gmail.com}

\nobibliography*

\usepackage[colorlinks=true]{hyperref}

\hypersetup{allcolors=[rgb]{0.1,0.1,0.4}}

\newtheorem{theorem}{Theorem} 

\newtheorem{conjecture}{Conjecture}
\newtheorem{properties}{Properties}
\newtheoremrm{rem}{Remark}

\mathrmdef{Hom}
\mathbfdef{Set}

\xymatrixrowsep{3pc}
\xymatrixcolsep{3pc}

\UseAllTwocells
\input diagxy

\SelectTips{eu}{}

\newcommand{\verticalize}[1]{\mathrel{\ooalign{\ensuremath{#1}\cr\hfil$\mid$\hfil}}}
\newcommand{\lorth}{\verticalize{\wedge}}
\newcommand{\oplorth}{\verticalize{\lor}}
\newcommand{\lordcell}[3]{\ar@{{}{ }{}}@/#1mm/[#2]^[#3]{\leq}}
\newcommand{\gordcell}[3]{\ar@{{}{ }{}}@/#1mm/[#2]^[#3]{\geq}}

\newcommand{\nat}[1]{\mathbb{#1}}
\newcommand{\cat}[1]{\mathscr{#1}}
\newcommand{\mor}[1]{\mathcal{#1}}
\newcommand{\san}[1]{\mathsf{#1}}

\newcommand{\xyc}[1]{\xymatrixcolsep{#1 pc}}

\newcommand{\ton}[1]{ \left( #1 \right)}
\newcommand{\qua}[1]{ \left[ #1 \right]}
\newcommand{\bra}[1]{ \left\{ #1 \right\}}

\newcommand{\id}[1]{\textup{id}_{#1}}
\newcommand{\Id}[1]{\textup{Id}_{#1}}
\newcommand{\harpup}{\ar@{-^{>}}}
\newcommand{\harpdown}{\ar@{-_{>}}}

\newcommand{\armo}{\ar@{>->}}
\newcommand{\arep}{\ar@{->>}}
\newcommand{\incl}{\ar@{^{(}->}@<-0.1pc>}
\newcommand{\emp}{\textup{\O}}
\newcommand{\prar}{\ar|-*=0@{|}}
\newcommand{\Epi}{\textup{\textbf{Epi}}}
\newcommand{\Mono}{\textup{\textbf{Mono}}}

\newcommand{\Clax}{\cat{C}^2_{\textup{lax}}}
\newcommand{\Coplax}{\cat{C}^2_{\textup{oplax}}}
\newcommand{\Ord}{\mathbf{Ord}}
\newcommand{\trat}{\textup{-}}
\newcommand{\Pcat}[1]{\mor{P}\ton{#1}}

\counterwithin*{equation}{section}

\begin{document}

\maketitle
\begin{abstract}
In this work we discuss a new type of factorisation systems for \textbf{Ord}-enriched categories. We start by defining the new notion of lax weak orthogonality, which involves the existence of lax diagonal morphisms for lax squares. Using the usual theory of factorisation systems as a blueprint, we introduce a lax version of weak and functorial factorisation systems. We provide a characterisation of lax functorial weak factorisation systems, namely lax weak factorisation systems such that their factorisations are lax functorial. Then, we present a particular case of such lax functorial weak factorisation systems which are equipped with additional lax monad structures. We finally explore some examples of these factorisation systems for categories of partial maps equipped with natural $\Ord$-enrichments.  We will first construct a particular lax algebraic weak factorisation system for categories of partial maps that isolates the total datum and the partial domain of a partial map. Then we will discuss the relation between factorisation systems on the base category and oplax factorisation systems on the induced category of partial maps.
\end{abstract}

\section*{Introduction}

Our study of lax and oplax factorisation systems stems from the introduction of a new
definition of orthogonality that encompasses diagonal liftings for lax or oplax squares, based on an Ord-enrichment of our category.

It differs substantially from the lax orthogonal factorisation systems studied in Ord-enriched categories in \cite{clementino2020lax,franco2019cofibrantly}, although in both cases the orthogonality condition does not impose uniqueness of the liftings. On one hand ours impose the existence of liftings for lax squares, and on the other hand theirs guarantee the existence of universal liftings.

In this paper we present a new type of factorisation system for \textbf{Ord}-categories. This study starts with the introduction of an orthogonality relation between morphisms that provides diagonal morphisms also for the broader class of lax squares. On the basis of this new definition we rebuild the already developed theory for ordinary weak factorisation systems introduced in \cite{freyd1972categories} and described in \cite{adamek2002weak,garner2009understanding,riehl2011algebraic}. Then, we analyse these classes and the equivalences that constitute their intersection.

We continue in the second section by presenting how these factorisations can be defined functorially providing the definition of lax functorial factorisation systems and pointing out the differences that arise from the ordinary case.

We follow on to the third section to describe how the former two factorisation systems interact. We provide a characterisation of those lax functorial weak factorisation systems. We also give a description of the two classes of morphisms of the underlying lax weak factorisation systems.

In the fourth section we present the lax case of the algebraic weak factorisation systems introduced in \cite{grandis2006natural,garner2009understanding} and we prove that in this setting they are lax weak and lax functorial as well.

Then we will focus on providing a broad description of the lax and oplax factorisation systems that arise for partial maps equipped with well-known definitions of order. The work on categories of partial maps and especially in their $\Ord$-enrichment is also inspired by the works of Fiore in \cite{fiore2004axiomatic}.

We will first construct a particular lax algebraic weak factorisation system for categories of partial maps that isolates the total datum of and the partial domain.

Then we will describe the relation between oplax factorisation systems on $\cat{C}$ and oplax factorisation systems on the correspondent category of partial maps. We will show that one induces the other and vice versa and that ordinary functorial weak factorisations transfer also their functorial properties to oplax factorisations on partial maps.

\section{Lax weak orthogonality and lax weak factorisation systems}

In this section we introduce a new notion of orthogonality which will be the base for the rest of the work. Let $\cat{C}$ be an $\Ord$-enriched category. We denote by $\Clax$ the category whose objects are morphisms in $\cat{C}$ and morphisms are squares $\ton{u,v}: f \longrightarrow g$ of the type
\begin{equation}\label{laxsquare}
\vcenter{
\xymatrix{
A \ar[r]^u \gordcell{_2.5}{dr}{@!45} \ar[d]_{f} & C \ar[d]^g \\
B \ar[r]_{v} & D.
}}
\end{equation}
We will refer to these squares as \textit{lax squares}. Then we introduce an orthogonality relation that provides lax diagonal morphisms also for these lax squares.

\begin{definition}\label{laxorthogonality}
Two morphisms in $\cat{C}^2$ are said to be \textbf{laxly weak orthogonal}, denoted by $f \lorth g$, if, for every lax square $\ton{u,v}: f \longrightarrow g$, there exists a morphism $d : cod(f) \longrightarrow dom(g)$, such that 
\begin{align}
    &\vcenter{\xymatrix{
    A \ar[r]^u \ar@{{}{ }{}}@/^3mm/[d]^(.35)[left]{\geq} \ar[d]_{f} & C \ar[d]^g \\
    B \ar[ur]|{d} \ar[r]_{v} & D \ar@{{}{ }{}}@/^3mm/[u]^(.35)[left]{\geq}
    }}
    &
    &\begin{cases}
    u \leq d \cdot f \\
    g \cdot d \leq v.
    \end{cases}
\end{align}
\end{definition}

We first observe that this constitutes a generalisation of weak orthogonality as described in \cite{adamek2002weak,riehl2011algebraic}. In fact, whenever the partial order in $\cat{C}$ is discrete, the two definitions coincide. Analogously, for any class of morphisms $\mor{H}$ we can define its lax weak orthogonal complements as the classes containing exactly those morphisms which are laxly weakly orthogonal to $\mor{H}$, to the left ($\mor{H}^{\lorth}$) or to the right ($^{\lorth}\mor{H}$). Lax weak orthogonal complements are easily shown to be closed under composition.

\begin{definition}
A \textit{lax weak prefactorisation system} is a pair $\ton{\mor{L}, \mor{R}}$ of classes of morphisms such that $\mor{R} = \mor{L}^{\lorth}$ and $\mor{L} = {}^{\lorth}\mor{R}$. Moreover if any morphism $f \in \Clax$ has an $\ton{\mor{L}, \mor{R}}$-factorisation
\begin{equation}\vcenter{
    \xymatrix{
    A \ar[rr]^f \ar[dr]_{\mor{L}\ni l_f} && B \\
    & W_f \ar[ur]_{r_f \in \mor{R}} &
    }}
\end{equation}
    then $\ton{\mor{L}, \mor{R}}$ is said to be a \textbf{lax weak factorisation system} (\textsc{lwfs}).
\end{definition}

We remark that, given a \textsc{lwfs} $\ton{\mor{L}, \mor{R}}$, for any lax square $\ton{u,v}: f \longrightarrow g$ there exists a morphism $\delta$ as in the diagram
\begin{equation}\vcenter{
    \xymatrix{
    A \gordcell{_2.5}{dr}{@!45} \ar[r]^u \ar[d]_{l_f} & C \ar[d]^{l_g} \\
    W_f \gordcell{_2.5}{dr}{@!45} \ar[d]_{r_f} \ar@{-->}[r]^{\delta} & W_g \ar[d]^{r_g} \\
    B \ar[r]_v & D.
    }}
\end{equation}
The arrow $\delta$ is obtained by the lax weak orthogonality relation $l_f \lorth r_g$.

The first natural step is to identify those morphisms which are laxly weak orthogonal with respect to any other morphism in the category. We recall that two morphisms $f:A\rightarrow B$ and $g:B \rightarrow A $ constitute an adjunction $f \dashv g$, if we have  $\id{A} \leq g \cdot f$ and $f \cdot g \leq \id{B} $.

\begin{proposition}\label{selforth}
For a morphism $f \in \Clax$, the following are equivalent:
\begin{enumerate}
    \item $f \lorth f$;
    \item $f$ is a left adjoint morphism to some $f^*$;
    \item $f \lorth \Clax$;
    \item $\Clax \lorth f$.
\end{enumerate}
\end{proposition}
\proof
1.$\Rightarrow$2. For a morphism that is orthogonal to itself we may consider the identity lax square (which is actually commutative) and by self orthogonality it must admit a lax diagonal filler. It follows that the diagonal morphism will be the right adjoint and the two 2-cells in the two triangles will define the adjunction sought.

2.$\Rightarrow$3. Considering a lax square $\ton{u,v}: f \longrightarrow g$ it is a mere calculation to prove that $u \cdot f^{*}$ is a lax diagonal filler for the lax square.

2.$\Rightarrow$4. Analogously for a lax square $\ton{u,v}: g \longrightarrow f$ we have that $f^{*} \cdot v$ is a lax diagonal filler.

3.$\Rightarrow$1. and 4.$\Rightarrow$1. are trivial.
\endproof

As a consequence of this result, we have that left adjoint morphisms $\textbf{LA}\ton{\cat{C}}$ belong to any lax weak orthogonal complement and they constitute the intersection between the two classes of morphisms of any lax weak prefactorisation system.

We point out that uniqueness of such lax diagonal liftings is not granted in general. In fact, given a morphism $f$ satisfying the conditions of Proposition \ref{selforth}, we have that in any lax square $\ton{u,v}:f \longrightarrow f$ both $d_u= u \cdot f^*$ and $d_v = f^* \cdot v$ are suitable lax diagonal morphisms.

As for ordinary factorisation systems, one could easily prove the following result.

\begin{proposition}\label{smallobjectlax}
Given a class of morphisms $\mor{H}\subseteq \cat{C}^2$, then $\mor{H} \subseteq {}^{\lorth}\ton{\mor{H}^{\lorth}}$ and $\mor{H} \subseteq \ton{^{\lorth}\mor{H}}^{\lorth}$. Moreover $\ton{^{\lorth}\ton{\mor{H}^{\lorth}}; \mor{H}^{\lorth}}$ and $\ton{^{\lorth}\mor{H}; \ton{^{\lorth}\mor{H}}^{\lorth}}$ are lax prefactorisation systems.
\end{proposition}

\section{Lax functorial factorisations}

A natural step forward is to study factorisation systems for lax arrow categories that are functorial. To do so, we consider the composition functor applied to lax arrow categories
\begin{equation}
    \xyc{3}
    \xymatrix{
    \Clax \times_{\mathscr{C}} \Clax \ar@<1.4ex>[r]|-{\pi_1} \ar@<0.1ex>[r]|-{\left( - \cdot -\right)} \ar@<-1.4ex>[r]|-{\pi_2}  & \Clax .
}
\end{equation}

We point out that objects in $\Clax \times_{\mathscr{C}} \Clax$ are pairs of composable morphisms and the arrows are triples of morphisms as
\begin{equation*}
    \xymatrix{
    A \gordcell{_2.5}{dr}{@!45} \ar[r]^a \ar[d]_{f} & A' \ar[d]^{f'} \\
    B \gordcell{_2.5}{dr}{@!45} \ar[d]_{g} \ar[r]^{b} & B' \ar[d]^{g'} \\
    C \ar[r]^c & C'.
    }
\end{equation*}

Then the following definition becomes a natural translation from the ordinary factorisation systems.

\begin{definition}\label{laxfuncotorialfactsys}
A \textbf{lax functorial factorisation system} is a functor $F: \Clax \longrightarrow \Clax \times_{\mathscr{C}} \Clax$, such that $\ton{\trat\cdot \trat}F =\Id{\Clax}$.
\end{definition}

A lax functorial factorisation system is then determined by a section of the composition functor applied to $\Clax$. When composing such functor $F$ with the projections $\pi_{1}$ and $\pi_{2}$ one obtains the usual functors $L,R: \Clax \longrightarrow \Clax$, and $K: \Clax \longrightarrow \mathscr{C}$.

A lax functorial factorisation system induces also the natural transformations $\eta : \Id{} \Rightarrow R$ and $\varepsilon : L \Rightarrow \Id{} $. Differently from ordinary factorisation systems, these transformations are not strict in general, but only oplax. In fact, considering $\eta$, for any lax square $\ton{u,v}: f \longrightarrow g$ we have that 
\begin{align*}
    &\vcenter{\xymatrix{
    A \gordcell{_2.5}{dr}{@!45} \ar[r]^u \ar[d]_{f} & C \ar[d]^{g} \ar[r]^{Lg} \ar@{{}{ }{}}@/_2.5mm/[dr]^{\eta_g} & Kg \ar[d]^{Rg} \\ B \ar[r]_v & D \ar[r]_{\id{D}} & D}
    }
    &\leq &
    &\vcenter{\xymatrix{
    A  \ar[d]_{f} \ar[r]^{Lf} \ar@{{}{ }{}}@/_2.5mm/[dr]^{\eta_f} & Kf \gordcell{_2.5}{dr}{@!45} \ar[d]^{Rf} \ar[r]^{K\ton{u,v}} & Kg \ar[d]^{Rg} \\
    B \ar[r]_{\id{B}} & B \ar[r]_v & D,}}
\end{align*}
since $Lg \cdot u \leq K\ton{u,v} \cdot Lf$ by definition of lax functorial factorisation. This amounts to have that $\eta_g \cdot \ton{u,v} \leq R \ton{u,v} \cdot \eta_f$. Similarly one can prove that $\varepsilon$ is an oplax natural transformation as well.

\section{Lax functorial weak factorisations}
Our goal now is to interlink the two concepts as it already happens in the traditional setting. An ordinary weak factorisation system $\ton{\mor{L}, \mor{R}}$ is said to underlie a functorial factorisation system $\ton{F,L,R,K}$ if every functorial factorisation is also an $\ton{\mor{L}, \mor{R}}$-factorisation.

We already know that, if an ordinary \textsc{wfs} $\ton{\mor{L}, \mor{R}}$ underlies the functorial factorisation system $\ton{F,L,R,K}$, then $\ton{\mor{L}, \mor{R}} =\ton{L\trat \textup{coalg} , R\trat \textup{alg}}$. This means that $\mor{L}$ contains those morphisms whose right component is a split epimorphism and $\mor{R}$ those morphisms whose left component is a split monomorphism (see for instance \cite{riehl2011algebraic}).
On the blueprint of this idea we want to investigate the conditions under which a lax functorial factorisation system has an underlying \textsc{lwfs} and to give a description for the latter.

We fix a lax functorial factorisation system with components $\ton{F,L,R,K}$. We consider the two classes 
\begin{align}\label{lfwfsdef}
    \mor{L}_F &=\bra{f \vert f \lorth Rf} & 
    \mor{R}_F &= \bra{f \vert Lf \lorth f}.
\end{align}

Similarly to ordinary functorial \textsc{wfs}, we consider a lax functorial factorisation system such that for every $f$, $Lf \in \mor{L}_F$ and $Rf \in \mor{R}_F$.
More precisely, we will consider lax factorisation systems such that $Lf \lorth RLf$ and $LRf \lorth Rf$ for every morphism $f$; we will call such lax functorial factorisation systems \textit{predistributive}. The reason for this name is that the assumption amounts to a certain distributivity of the lax functorial factorisation system as depicted in \eqref{predistributivediag}.

\begin{proposition}\label{predistrprop}
Let $\ton{F,L,R,K}$ be a lax functorial factorisation system and $f$ any morphism. If $Rf \in \mor{R}_F$, then $f \in \mor{L}_F$ if and only if there exists a lax diagonal lifting $\rho_f$ in the square $\eta_f$. If $Lf \in \mor{L}_F$, then $f \in \mor{R}_F$ if and only if there exists a lax diagonal lifting $\lambda_f$ in the square $\varepsilon_f$. For any $f$, we will denote such lax liftings by
\begin{align}
    &\vcenter{
    \xymatrix{
        A \ar[r]^{Lf} \ar@{{}{ }{}}@/^3mm/[d]^(.35)[left]{\geq} \ar[d]_{f} & Kf \ar[d]^{Rf} \\
        B \ar[ur]|{\rho_{f}} \ar[r]_{\id{B}} & B \ar@{{}{ }{}}@/^3mm/[u]^(.35)[left]{\geq},
        }
    }
    &
    &\vcenter{
        \xymatrix{
        A \ar[r]^{\id{A}} \ar@{{}{ }{}}@/^3mm/[d]^(.35)[left]{\geq} \ar[d]_{Lf} & A \ar[d]^{f} \\
        Kf \ar[ur]|{\lambda_{f}} \ar[r]_{Rf} & B \ar@{{}{ }{}}@/^3mm/[u]^(.35)[left]{\geq}.
        }
    }
\end{align}
\end{proposition}
\proof
We prove only the first statement, since the second one follows by duality.

One direction is trivial since the existence of such a $\rho_{f}$ is a direct consequence of $f \lorth Rf$.

For the non-trivial implication, we need to prove that $f \lorth Rf$. We consider a lax square 
\begin{align}\label{Krectangle}
    &\vcenter{
    \xymatrix{
    A \ar[r]^u \gordcell{_2.5}{dr}{@!45} \ar[d]_{f} & Kf \ar[d]^{Rf} \\
    B \ar[r]_{v} & B
    }}&
    &\longmapsto & 
    &\vcenter{\xymatrix{
    A \gordcell{_2.5}{dr}{@!45} \ar[r]^u \ar[d]_{Lf} & Kf \ar[d]|{LRf} \\
    Kf \gordcell{_2.5}{dr}{@!45} \ar[d]|{Rf}  \ar[r]^{K\ton{u,v}} & Kg \ar[d]^{RRf} \ar@/_8mm/[u]_{\lambda_{Rf}} \\
    B \ar@/^8mm/[u]^{\rho_{f}} \ar[r]_v & B,
    }}
\end{align}
where $\rho_{f}$ is a diagonal morphism of $\eta_f$ existing by assumption and $\lambda_{Rf}$ is a diagonal morphism of $\varepsilon_{Rf}$ existing since $LRf \lorth Rf$.

We consider $\Delta=\lambda_{Rf} \cdot K \ton{u,v} \cdot \rho_{f}$.
Keeping in mind the definitions of $\rho_f$ and $\lambda_{Rf}$ as lax diagonal fillers of $\eta_f$ and $\varepsilon_{Rf}$, we have that
\begin{equation}
    u \leq \lambda_{Rf} \cdot LRf \cdot u \leq \lambda_{Rf} \cdot K \ton{u,v} \cdot Lf \leq \lambda_{Rf} \cdot K \ton{u,v} \cdot \rho_{f} \cdot Rf \cdot Lf = \Delta \cdot f
\end{equation}
and
\begin{equation}
    Rf \cdot \Delta = RRf  \cdot LRf \cdot \lambda_{Rf} \cdot K \ton{u,v} \cdot \rho_{f} \leq RRf  \cdot K \ton{u,v} \cdot \rho_{f} \leq v \cdot Rf \cdot \rho_{f} \leq v.
\end{equation}
This yields that $\Delta$ is the diagonal morphism sought and $f \lorth Rf$.
\endproof

\begin{remark}
We recall that, for a lax pointed functor $\ton{T, \theta}$, a pair $\ton{A, a}$ is said to be a \textit{lax algebra} if $\id{A} \leq a \cdot \theta_A$.
Then we point out that, for any morphism $f$, the existence of a lax diagonal filler for $\varepsilon_f$ is equivalent to have a lax algebra structure for the lax pointed endofunctor $\ton{R,\eta}$. In fact the lax square $\alpha = \ton{\lambda_f , \id{B}} : Rf \longrightarrow f$ is such that $\id{f} \leq \alpha \cdot \eta_f$. Similarly any morphism $g$ admits a lax diagonal morphism for $\varepsilon_g$ if and only if it has a lax coalgebra structure for $\ton{L, \varepsilon}$.
\end{remark}

\begin{corollary}
Let $\ton{F,L,R,K}$ be a predistributive lax functorial factorisation system. Then 
\begin{align}
\begin{split}
    \mor{L}_F & = \bra{f \vert \textup{ $\eta_f$ has a lax diagonal morphism}} = L \textup{-coalg}_{\textup{lax}} \\
    \mor{R}_F & = \bra{f \vert \textup{ $\varepsilon_f$ has a lax diagonal morphism}}= R \textup{-alg}_{\textup{lax}}.
\end{split}
\end{align}
\end{corollary}

We observe that, given a morphism $f$ that lies in $\mor{L}_F \cap \mor{R}_F$, then the composition $\lambda_f \cdot \rho_f$ is a right adjoint to $f$.

According to these remarks, the choice of the name predistributive points to the existence, for any morphism $f$, of the following diagram
\begin{equation}\label{predistributivediag}
    \vcenter{\xymatrix{
    Kf \ar@{{}{ }{}}@/^3mm/[d]^(.6)[left]{\geq} \ar@{=}[dr] \ar[d]_{LRf} \ar[r]^{\rho_{Lf}} & KLf \ar[d]^{RLf} \\
    KRf \ar[r]_{\lambda _{Rf}} & Kf \ar@{{}{ }{}}@/^3mm/[u]^(.6)[left]{\geq} 
    }}
\end{equation}
and it coincides with the assumption that for every morphism $f$, $Rf$ is a lax $\ton{R , \eta}$-algebra and $Lf$ is a lax $\ton{L, \varepsilon}$-coalgebra.

This diagram resembles the distributivity transformation described in the next section, even if it carries less structure. In fact, it is not in general a natural transformation and it does not satisfy any distributivity law as described in \cite{bourke2016algebraic,clementino2016lax} or in \cite[Section 4]{clementino2020lax}.

\begin{theorem}\label{theorem}
Let $\ton{F,L,R,K}$ be a predistributive lax functorial factorisation system. Then $\ton{\mor{L}_F, \mor{R}_F}$ is a lax weak factorisation system. Moreover, for any lax functorial weak factorisation system $\ton{\mor{L},\mor{R}}$ with lax functorial factorisation $\ton{F,L,R,K}$, $\ton{\mor{L},\mor{R}} = \ton{\mor{L}_F, \mor{R}_F}$.
\end{theorem}
\begin{proof}
We start by proving that $\mor{L} \lorth \mor{R}$. Let $f \in \mor{L}$ and $g \in \mor{R}$. We have that $f \lorth Rf$ and $Lg \lorth g$, by the existence of the two morphisms $\rho_{f}$ and $\lambda_g$.
We factorise a lax square as
\begin{multline}\label{theorectangle}
    \vcenter{
    \xymatrix{
    A \ar[r]^u \gordcell{_2.5}{dr}{@!45} \ar[d]_{f} & C \ar[d]^g \\
    B \ar[r]_{v} & D
    }}
    \longmapsto 
    \vcenter{\xymatrix{
    A \gordcell{_2.5}{dr}{@!45} \ar[r]^u \ar[d]_{Lf} & C \ar[d]^{Lg} \\
    Kf \gordcell{_2.5}{dr}{@!45} \ar[d]_{Rf} \ar[r]^{K\ton{u,v}} & Kg \ar[d]^{Rg} \ar@/_8mm/[u]_{\lambda_{f}} \\
    B \ar@/^8mm/[u]^{\rho_{f}} \ar[r]_v & D.
    }}
\end{multline}
The morphism $\Delta= \lambda_g \cdot K\ton{u, v} \cdot \rho_f$ is a lax diagonal morphism for the lax square taken into account.
In conclusion this yields that $\mor{L} \lorth \mor{R}$.

Moreover, for any $f \lorth \mor{R}$, it follows that $f \lorth Rf$, since $Rf \in \mor{R}$ by lax predistributivity, which induces that $f \in \mor{L}$, namely $\mor{R}^{\lorth} \subseteq \mor{L}$. By an analogous argument $\mor{L}^{\lorth} \subseteq  \mor{R}$ holds. 

For the second claim we consider $\ton{\mor{L}, \mor{R}}$ a functorial weak factorisation system. We have have that
\begin{align*}
    f & \in \mor{L} & &\Leftrightarrow & f & \in ^{\lorth}\mor{R} & &\Rightarrow & &f \lorth Rf & &\Leftrightarrow & f &\in \mor{L}_F \\
    f & \in \mor{R} & &\Leftrightarrow & f & \in \mor{L}^{\lorth} & &\Rightarrow & & Lf \lorth f & &\Leftrightarrow & f &\in \mor{R}_F
\end{align*}
so $\mor{L} \subseteq \mor{L}_F$ and $\mor{R} \subseteq \mor{R}_F$ and since inclusion is dual for lax weak orthogonal classes, this implies that $\ton{\mor{L},\mor{R}} = \ton{\mor{L}_F, \mor{R}_F}$.
\end{proof}

In conclusion, the last theorem gives us a description of lax \textsc{fwfs} as algebras and coalgebras of the functorial components. Moreover, it states that if every $L$-component is a lax coalgebra and every $R$-component is a lax algebra, then they form the unique \textsc{lwfs} underlying the given \textsc{fwfs}.

\section{Lax algebraic weak factorisations}
In this section we will present a class of functorial factorisation systems which satisfy the condition above and come equipped with a richer structure close to the one of a monad. This construction parallels in this lax context that of algebraic weak factorisation systems as in \cite{grandis2006natural,garner2009understanding,bourke2016algebraic}.

We recall the definition of lax monad that we reprise from \cite{bunge1974coherent}. We define here an $\Ord$ version of this definition and, although this work refers to lax natural transformations, we remark that it is actually the same type of transformation we call oplax according to what appears to be the most used choice in literature. The only difference is that we will use a definition that involves usual functors and not lax functors, since it is the particularisation that best fits our purposes.

\begin{definition}
For an $\Ord$-enriched category $\cat{C}$, a lax monad is a triple $\ton{T, \eta , \mu }$, such that 
\begin{itemize}
    \item $T: \cat{C} \longrightarrow \cat{C}$ \textup{ is a functor;}
    \item $\eta : \Id{} \Longrightarrow T$ \textup{ is an oplax natural transformation;}
    \item $\mu : RR \Longrightarrow R $ \textup{is an oplax natural transformation;}
\end{itemize}
and such that the following lax monads laws are satisfied
\begin{align}\label{laxmonadlaws}
    &\vcenter{
    \xymatrix{
    T \ar[dr]_{\id{T}} \ar[r]^{T\eta} \lordcell{^2}{dr}{@!45} & TT \ar[d]|{\mu} & T \ar[ld]^{\id{T}} \ar[l]_{\eta_T} \lordcell{_8}{dl}{@!-45} \\
    &T &
    }
    }
    &
    &\vcenter{
    \xymatrix{
    TTT \ar[d]_{\mu_T} \ar[r]^{T\mu} \gordcell{_2.5}{dr}{@!45} & TT \ar[d]^{\mu} \\
    TT \ar[r]_{\mu}  & T.
    }
    }
\end{align}
\end{definition}

We now define the factorisation systems we are interested in.

\begin{definition}
A \textup{lax algebraic weak factorisation system} (\textsc{lawfs}) is a functorial factorisation system $\ton{F,L,R,K}$ such that $\ton{R,\eta}$ is part of a lax monad $\ton{R,\eta , \Theta}$, $\ton{L,\varepsilon}$ is part of a lax comonad $\ton{L,\varepsilon , \Omega }$, and there exists a distributivity law $\Delta : LR \Longrightarrow RL$ of the comonad over the monad in the sense that the following diagram commutes
\begin{equation}\label{distriblawdiag}
    \vcenter{
    \xymatrix{
    LRR \ar[d]_{L\Theta} \ar[r]^{\Delta_R} & RLR \ar[r]^{R \Delta} & RRL \ar[d]^{\Theta_L} \\
    LR \ar[rr]^{\Delta} \ar[d]_{\Omega_R} && RL \ar[d]^{R\Omega} \\
    LLR \ar[r]^{L\Delta} & LRL \ar[r]^{\Delta_L} & RLL.
    }
    }
\end{equation}
\end{definition}

As said before, these factorisation systems constitute a subclass of lax functorial weak factorisation systems as we prove in the following proposition.

\begin{proposition}
A \textsc{lawfs} $\ton{F,L,R,K}$ is a lax functorial weak factorisation system.
\end{proposition}
\proof
We can prove the statement by showing that $\ton{F,L,R,K}$ is lax predistributive. Let $f$ be any morphism. Then we consider the lax square given by $\theta_f$, whose defining 2-cell is $Rf \cdot \theta_{f} \leq RRf$. Then, by the lax monad law $\id{R} \leq \Theta \cdot \eta_R$ in \eqref{laxmonadlaws}, we can deduce, restricting it to the domains, that $\id{Kf} \leq \theta_{f} \cdot LRf$. Thus we have that $\theta_{f}$ is a lax diagonal morphism for $\varepsilon_{Rf}$.
The same argument on the comonad yields that $\omega_f$, the codomain morphism of the comultiplication $\Omega$ of the comonad, is the lax diagonal morphism for $\eta_{Lf}$. Now we want to check that $Lf \lorth RLf$ and $LRf \lorth Rf$. We will prove only the first one, since the arguments for the second are similar.
We consider any lax square between $Lf $ and $ RLf$ and its factorisation

\begin{multline}\label{therectangle}
    \vcenter{
    \xymatrix{
    A \ar[r]^u \gordcell{_2.5}{dr}{@!45} \ar[d]_{Lf} & KLf \ar[d]^{RLf} \\
    Kf \ar[r]_{v} & Kf
    }}
    \longmapsto 
    \vcenter{\xymatrix{
    A \gordcell{_2.5}{dr}{@!45} \ar[r]^u \ar[d]_{LLf} & KLf \ar[d]|{LRLf} \\
    KLf \gordcell{_2.5}{dr}{@!45} \ar[d]|{RLf} \ar[r]^{K\ton{u,v}} & KRLf \ar[d]^{RRLf} \ar@/_8mm/[u]_{\theta_{Lf}} \\
    Kf \ar@/^8mm/[u]^{\omega_{f}} \ar[r]_v & Kf.
    }}
\end{multline}

Then $\theta_{Lf} \cdot K \ton{u,v} \cdot  \omega_f$ is a lax diagonal morphism. In fact, due to the rules of the monad, we have
\begin{equation}
  \begin{cases}
u \leq \theta_{Lf} \cdot LRLf \cdot u \leq \theta_{Lf} \cdot K \ton{u,v} \cdot LLf \leq \theta_{Lf} \cdot K \ton{u,v} \cdot \omega_f \cdot RLf \cdot LLf   \\
RRf \cdot LRLf \cdot \theta_{Lf} \cdot K \ton{u,v} \cdot \omega_f \leq RRLf \cdot K\ton{u,v} \cdot \omega_f \leq v \cdot RLf \cdot \omega_f \leq v.
\end{cases}  
\end{equation}
This implies that, for every $f$, we have that $Lf \lorth RLf$, and similarly $LRf \lorth Rf$. Therefore $\ton{F,L,R,K}$ is lax predistributive and hence it is a lax functorial weak factorisation system.
\endproof

We remark that in general a lax predistributive functorial factorisation system does not yield a complete distributivity law, since it is not even true that the square \eqref{predistributivediag} is a natural transformation, and moreover we do not have the existence of the 2-cells in the distributivity law \eqref{distriblawdiag}.

\begin{remark}
We observe that what we have described for $\Clax$ can be expressed in a dual fashion for $\Coplax$. We can then define orthogonality for oplax squares and oplax factorisation systems. We remark that the equivalences for oplax weak orthogonality are right adjoint morphisms and that, for an oplax functorial factorisation system, $\eta$ and $\varepsilon$ are lax natural transformations and the monads involved in oplax \textsc{awfs} will be obviously oplax monads. We have then an equally powerful dual set of results that can be used for $\Coplax$. We will denote the oplax weak orthogonality relation by $\oplorth$.
\end{remark}

\section{Categories of partial maps}
This second part of the work is dedicated to the study of lax and oplax factorisation systems for categories of partial maps. First we recall in this section some useful definitions, notations and properties of categories of partial maps. Some of the features that we present may be found in \cite{fiore2004axiomatic,robinson1988categories}. 

Let $\cat{C}$ be a category  with a good notion of subobjects, i.e. a class of monomorphisms $\mor{S}$ that is closed under composition, pullback stable, and that contains all sections (\textsc{ri}). We consider $\cat{C}$, with pullbacks along morphisms in $\mor{S}$. Then a partial map is a span

\begin{equation*}
    \xymatrix{
    D_f \ar@{>->}[d]_{\sigma_f} \ar[dr]^{\varphi_f} & \\
    A \prar[r]_f & B
    }
\end{equation*}
with $\sigma_f$ in $\mor{S}$ and $\varphi_f$ a generic morphism in $\cat{C}$. Composition among partial maps operates via pullback. We will often use hereon the notation $D_{\trat}$, $\sigma_\trat$ and $\varphi_\trat$ to refer to the partial domain and the partial components of a partial map.

The category obtained will be denoted by $\mor{P}_{\mor{S}}\ton{\cat{C}}$ or $\Pcat{\cat{C}}$, where it cannot generate any ambiguity.
Categories of partial maps have a partial order between morphisms defined by $f\preceq g$ if $f$ is a domain restriction of $g$; formally if there exists $s \in \mor{S}$ making the following diagram commute:
\begin{equation}\label{partialmaporder}
    \vcenter{
    \xymatrix{
    D_f \ar@{>->}[d]_{s} \ar@/_2pc/[dd]_{\sigma_f} \ar[ddr]^{\varphi_f} \ar@{{}{ }{}}@/_1mm/[ddr]_{\blacktriangle} & \\
    D_g \ar@{>->}[d]_{\sigma_g} \ar[dr]|{\varphi_g} & \\
    A \prar@<1.5mm>[r]^f \prar@<-1.5mm>[r]_g & B.
    }}
\end{equation}

Furthermore if there exists a partial order $\sqsubseteq $ among maps in $\cat{C}$, then we can define a partial order on partial maps induced by the one on the base category. In fact, if $\cat{C}$ is $\Ord$-enriched and all maps in $\mor{S}$ are full (a map $f$ is full if any 2-cell $f \cdot u \sqsubseteq f \cdot v$ yields the 2-cell $u \sqsubseteq v$), then we have another order relation $f \leq g$ if and only if $\blacktriangle$ in \eqref{partialmaporder} is the 2-cell $\varphi_f \sqsubseteq \varphi_g  \cdot  s$. It is clear that if we consider a discrete $\Ord$-enrichment on $\cat{C}$, then the two partial orders coincide.

\begin{notation}
We will distinguish the two order relations on a category $\Pcat{\cat{C}}$ with the symbols used above. Thus $\preceq$ will denote the order induced by equalities, while $\leq$ will denote an order induced by any $\Ord$-enrichment on $\cat{C}$.
\end{notation}

\begin{notation}
Given the 2-cells $f \leq g$ and $f' \leq g'$ properly composable, we will denote by $o' * o$ the subobject morphism defining the composition 2-cell $f' \cdot f \leq g' \cdot g $. Moreover, we will denote by $v * o$ the subobject morphism defining the composition 2-cell with an identity 2-cell such as $v \cdot f \leq v \cdot g$.
\end{notation}

We will now state some facts that we will frequently evoke in the later discussion. We recall that a partial morphism $f$ is said to be \textit{total} if its partial domain $D_f$ is the whole domain of $f$.

\begin{properties}\label{orderproperties}
Given a partial order on $\Pcat{\cat{C}}$, we have the following properties:
\begin{itemize}
    \item if $f$ is total and $f \leq g$, then $g$ is total. Moreover, if the partial order is induced by the discrete one, then $f=g$. Hence total maps are maximal elements in their own Hom-Sets;
    \item for any composition of partial morphisms $g  \cdot  f$, it is true that $D_{g  \cdot  f} \rightarrowtail D_f$. Furthermore, if $f  \cdot  g$ is total, then $g$ is total.
\end{itemize}
\end{properties}

\subsection{Adjunctions between partial maps}
Since adjoint morphisms constitute the lax equivalences related to lax factorisation systems, we proceed giving a description of adjoint morphisms for categories of partial maps.

\begin{proposition}\label{adjprop}
Let $\Pcat{\cat{C}}$ be a category of partial maps with $\cat{C}$ an $\Ord$-enriched category. A pair of morphisms in $\Pcat{\cat{C}}$ constitute an adjunction $f \dashv g$ if and only if $f$ is total and $\varphi_f = \sigma_g  \cdot  \Tilde{\varphi_f}$ such that $\Tilde{\varphi_f} \dashv \varphi_g$ in $\cat{C}$.
\end{proposition}
\proof
First we consider the adjunction $f \dashv g$ in $\Pcat{\cat{C}}$. By the properties in \ref{orderproperties}, the 2-cell $\id{A} \leq g  \cdot  f$ yields that $g  \cdot  f$ is total and therefore $f$ is total. We write explicitly the 2-cells
\begin{align}\label{explicitadjunc}
    &\vcenter{
    \xymatrix{
    A \ar@{>->}[d] \ar@/_2pc/[dd]|{\id{A}} \ar[dr]|{\Tilde{\varphi_f}} \ar@/^3.4pc/[ddrr]|{\id{A}} \ar@{{}{ }{}}@/^1pc/[ddrr]^[@!45]{\sqsupseteq} && \\
    A \ar@{>->}[d] \ar[dr]|{\varphi_f} & D_g \ar@{>->}[d]|{\sigma_g} \ar[dr]|{\varphi_g} & \\
    A \prar[r]_f & B  \prar[r]_g & A.
    }
    }
    &
    &\vcenter{
    \xymatrix@R=1pc@C=1.5pc{
    B  \ar@/^4pc/[dddddrrrrr]|{\id{B}} \ar@/_2pc/[dddddr]|{\id{B}}  \ar@{{}{ }{}}@/^1.2pc/[dddddrrrrr]^[@!45]{\sqsubseteq} && \\
    & D_g \armo[ul]|{\sigma_g} \ar@{>->}[dd] \ar[ddrr]|{\varphi_g}  &&&& \\
    &&&&& \\
    & D_g \ar@{>->}[dd]|{\sigma_g} \ar[ddrr]|{\varphi_g} && A \ar@{>->}[dd] \ar[ddrr]|{\varphi_f} && \\
    \\
    & B \prar[rr]_g && A  \prar[rr]_f && B.
    }
    }
    \end{align}
The subobject morphism of the left diagram is $\id{A}$ due to the stated totality of $f$ and $g  \cdot  f$. We notice immediately that $\varphi_f = \sigma_g  \cdot  \Tilde{\varphi_f}$, where $\Tilde{\varphi_f}$ is the pullback of $\varphi_f$ along $\sigma_g$. Furthermore, we have the 2-cells $\id{A} \sqsubseteq \varphi_g  \cdot  \Tilde{\varphi_f}$ and $\sigma_g  \cdot  \Tilde{\varphi_f}  \cdot  \varphi_g = \varphi_f  \cdot  \varphi_g \sqsubseteq \sigma_g$. Since by assumption $\sigma_g \in \mor{S}$ is full, we deduce that $\Tilde{\varphi_f}  \cdot  \varphi_g \sqsubseteq \id{B}$. This yields the existence of the adjunction $\Tilde{\varphi_f} \dashv \varphi_g$ in $\cat{C}$.
The other direction is proved simply by calculation. In fact, due to the rules of composition among partial maps, one obtains explicitly the diagrams in \eqref{explicitadjunc}, henceforth proving $f \dashv g$ in $\Pcat{\cat{C}}$.
\endproof

If we consider an ordinary category $\cat{C}$, we can state the following corollary.

\begin{corollary}
Let $\Pcat{\cat{C}}$ be a category of partial morphisms equipped with the $\Ord$-enrichment induced by the ordinary category $\cat{C}$. A pair of morphisms in $\Pcat{\cat{C}}$ constitute an adjunction $f \dashv g$ if and only if $f=\ton{\id{A} , \sigma}$ and $g=\ton{\sigma , \id{A}}$, for some $\sigma \in \mor{S}$.
\end{corollary}

This is due to the fact that pairs of adjoint morphisms in a discrete category are pairs of inverse morphisms.

We conclude this discussion on adjoint partial maps with a remark on the relation between the adjunctions of partial maps and the adjunctions in the base $\Ord$-category. First we can define a partial order for adjoint pairs in $\Pcat{\cat{C}}$ such that $\ton{f \dashv g} \leq \ton{f' \dashv g'}$ if and only if there exists $s \in \mor{S}$ such that $\varphi_f = s  \cdot  \varphi_{f'}$ and $\sigma_{g}= s  \cdot  \sigma_{g'}$.
It is straightforward to check that any pair of adjoint maps in $\cat{C}$ still constitutes an adjoint pair of total maps in $\Pcat{\cat{C}}$. Vice versa, Proposition \ref{adjprop} induces a process to obtain an adjoint pair in $\cat{C}$ from any adjoint pair in $\Pcat{\cat{C}}$ and in particular it preserves adjunction on total maps. Following the notation of the Proposition it becomes easy to check that $\ton{f \dashv g} \leq \ton{\ton{\id{A} , \Tilde{\varphi_f}} \dashv \ton{\id{B} , \varphi_g}}$. It yields an adjunction between the partially ordered sets $\textbf{Adj} \ton{\Pcat{\cat{C}}}$ and $\textbf{Adj} \ton{\cat{C}}$.

\section{Domain-Total factorisation}\label{TotDom}
In this section we will provide a general construction of a lax functorial weak factorisation system for any category of partial maps equipped with the most general definition of partial order described above. We consider and $\Ord$-category $\cat{C}$ and the category of partial maps $\Pcat{\cat{C}}$ induced by any appropriate class of subobject $\mor{S}$. Given any partial map $f$ in $\Pcat{\cat{C}}$ we can factorise it as
\begin{align}\label{totinjfact}
    &\vcenter{\xymatrix{
    D_f \armo[d]_{\sigma_f} \ar[dr]^{\varphi_f} & \\
    A \prar[r]_f & B
    }}
    &
    &=
    &
    &\vcenter{\xymatrix{
    D_f \armo[d]_{\sigma_f} \ar[dr]|{\id{D_f}} &  D_f \armo[d]|{\id{D_f}} \ar[dr]^{\varphi_f} & \\
    A \prar[r]_{Lf} & D_f \prar[r]_{Rf} & B.
    }}
\end{align}
Our goal now will be to prove that such factorisation is a lax functorial weak factorisation system. The first step is proving that it is a lax functorial factorisation, namely that $L$ and $R$ are part of a functorial factorisation $F : \Pcat{\cat{C}}^2_{\text{lax}} \longrightarrow \Pcat{\cat{C}}^2_{\text{lax}} \times_{\Pcat{\cat{C}}} \Pcat{\cat{C}}^2_{\text{lax}}$. We consider a lax square $\ton{u,v}: f \longrightarrow g$ and its factorisation through $L $ and $ R$

\begin{align}\label{functdiagrams}
    &\vcenter{\xymatrix{
    A \prar[r]^u \gordcell{_2.5}{dr}{@!45} \prar[d]_{f} & C \prar[d]^g \\
    B \prar[r]_{v} & D
    }}
    &
    & \longmapsto
    &
    &\vcenter{\xymatrix{
    A \ar@{{}{ }{}}@/_2.5mm/[dr]^{\fbox{1}} \ar[r]^u \ar[d]_{Lf} & C \ar[d]^{Lg} \\
    D_f \ar@{{}{ }{}}@/_2.5mm/[dr]^{\fbox{2}} \ar[d]_{Rf} \ar@{-->}[r]^{K\ton{u,v}} & D_g \ar[d]^{Rg} \\
    B \ar[r]^v & D
    }}
    &
    & \text{where }K\ton{u,v} \text{ is}
    &
    &\vcenter{
    \xymatrix@R=1.5pc{
    D_{g  \cdot  u} \armo[d] \armo[ddr]^{\overline{\varphi_u} = \sigma_{g}^{*}\ton{\varphi_u}} & \\
    D_{v  \cdot  f} \armo[d] & \\
    D_f \prar[r] & D_g.
    }}
\end{align}

In the second diagram the upper square \fbox{1} is actually commutative, while it is a simple calculation that the subobject morphism $o:D_{g  \cdot  u} \rightarrowtail D_{v  \cdot  f}$ defining $g  \cdot  u \leq v  \cdot  f$ is also a witness that square \fbox{2} is a lax square. We consider two composable lax squares $\ton{u , v}: f \longrightarrow g$ and $\ton{u' , v'}: g \longrightarrow h$. We have to prove that $K \ton{u'  \cdot  u , v'  \cdot  v} = K\ton{u', v'}  \cdot  K \ton{u,v}$. To prove this equality one has to apply the composition rules among partial maps while remembering the definition of $K$ in \eqref{functdiagrams}. Then the key argument is that the two arrows 
\begin{equation*}
    \vcenter{
    \xymatrix@R=1.5pc@C=1.5pc{
    D_{h  \cdot  u'} \armo[rr]^{o}\armo[d] && D_{v'  \cdot  g} \armo[d] \\
    D_{u'} \armo[dr] && D_{g} \armo[dl] \\
    & A'
    }
    }
\end{equation*}
depict the same subobject equivalence class of $D_{h  \cdot  u'}$ in $A'$ and therefore pullbacks of $\varphi_u$ along these arrows are equal. This yields our thesis and shows that $L$ and $R$ constitute a lax functorial factorisation system.

Next we discuss lax predistributivity of $\ton{F , L, R, K}$. We proceed giving in our setting a particular description of the two classes.

First we show that $\mor{L}_F$ is the class $\textbf{PLA}$ of partial left adjoint morphisms, namely any morphism $f$ such that $\varphi_f= \sigma  \cdot  \overline{\varphi_f}$, $\sigma \in \mor{S}$ and $\Tilde{\varphi_f}$ is a left adjoint in $\cat{C}$. It is clear that if $f \in \textbf{PLA}$, then $Rf$ is a left adjoint morphism and therefore $f \lorth Rf$. On the other hand if $f \in \mor{L}_F$, then we have a lax diagonal morphism $\delta$ for the square denoted by $\eta_f$. Writing the 2-cells explicitly one can see easily that $\varphi_f = \sigma_{\delta}  \cdot  \overline{\varphi_f}$ and that $\overline{\varphi_f} \dashv \varphi_{\delta}$. Hence we have that $\mor{L}_F = \textbf{PLA}$.

As for $\mor{R}_F$ we notice first that it contains all total maps. In fact, if $f$ is total, then $Lf$ is an identity, therefore $Lf \lorth f$. If we consider $f$ in $\mor{R}_F$, again we have a lax diagonal morphism $\delta '$ for the square denoted by $\varepsilon_f$, which in particular yields that $\id{A} \leq \delta  \cdot  Lf$. Hence $Lf$ is total. This implies that $\mor{R}_F = \textbf{Tot}$.

Using the new characterizations of $\mor{L}_F$ and $\mor{R}_F$ given above, to prove lax predistributivity it is enough to notice that, for every partial map $f$, $Lf \in PLA$ and $Rf \in \textbf{Tot}$, which is trivially true. This yields that $\ton{\textbf{PLA}, \textbf{Tot}}$ is a lax weak factorisation system underlying the lax functorial factorisation system above.

We conclude remarking that, if $\cat{C}$ is a discrete category, then $\textbf{PLA}$ contains those morphisms such that $\varphi_f \in \mor{S}$. This is due to the fact that the only adjoint pairs are isomorphisms. We denote such class of morphisms by $\overline{\mor{S}}$.

\section{Oplax factorisation systems on partial maps and ordinary factorisation systems on total maps}
In the following section we will display the close links between factorisation systems on a category $\cat{C}$ and oplax factorisation systems. In the first paragraph we will focus on describing a procedure that produces an oplax weak factorisation system on $\Pcat{\cat{C}}$ from a stable oplax weak factorisation system on $\cat{C}$. Then we will analyse how functoriality is transferred to such factorisation systems. Thereafter we will proceed to study how oplax weak factorisation systems on partial maps may be restricted to factorisation systems among total maps.

\subsection{Oplax \textsc{wfs}s from total maps to partial maps}
We consider an $\Ord$-category $\cat{C}$ equipped with an oplax \textsc{wfs} $\ton{\mor{E}, \mor{M}}$ such that $\mor{E}$ is a class of morphisms stable under pullbacks along morphisms in $\mor{S}$. Looking at $\Pcat{\cat{C}}$ we can factorise each partial map $f$ as
\begin{align} \label{EMfactdiag}
    &\vcenter{\xymatrix{
    D_f \armo[d]_{\sigma_f} \ar[dr]^{\varphi_f} \ar[r]^{e_{\varphi_f}} &  M_f \ar[d]^{m_{\varphi_f}} \\
    A \prar[r]_f & B
    }}
    &
    &=
    &
    &\vcenter{\xymatrix{
    D_f \armo[d]_{\sigma_f} \ar[dr]|{e_{\varphi_f}} &  M_f \armo[d]|{\id{M_f}} \ar[dr]^{m_{\varphi_f}} & \\
    A \prar[r]_{e_f} & M_f \prar[r]_{m_f} & B.
    }}
\end{align}
We consider the following classes of partial morphisms
\begin{align}
    \overline{\mor{E}}&=\bra{f \vert \varphi_f \in \mor{E}} &
    \overline{\mor{M}}&=\bra{f \vert \varphi_f \in \mor{M}}.
\end{align}
We will prove that it constitutes an oplax weak factorisation system. First we prove that $\overline{\mor{E}} \oplorth \overline{\mor{M}}$. We consider $f \in \overline{\mor{E}} $, $g \in \overline{\mor{M}}$ and the oplax square $\ton{u,v}: f \longrightarrow g$. Writing explicitly the oplax square, we have that $\sigma_{v}^*\ton{\varphi_f} $ is in $ \mor{E}$ by the condition of stability under pullbacks. By oplax weak orthogonality of $\ton{\mor{E}, \mor{M}}$ there exists an oplax diagonal filler $d$ for the oplax square
\begin{equation}\label{stablediag}
    \xymatrix{
    D_{v  \cdot  f} \ar@{{}{ }{}}@/^8mm/[d]^(.40)[left]{\sqsubseteq } \ar[d]_{\mor{E} \ni \sigma_{v}^*\ton{\varphi_f}} \armo[r]^{o} & D_{g  \cdot  u} \ar[r]^{\sigma_{g}^*\ton{\varphi_u}} & D_g \armo[d]^{\varphi_g \in \mor{M}} \\
    D_v \ar@{-->}[urr]|{\exists d} \ar[rr]_{\varphi_v} && D, \ar@{{}{ }{}}@/^8mm/[u]^(.40)[left]{\sqsubseteq}
    }
\end{equation}
where $o$ is the morphism in $\mor{S}$ involved in the definition of the 2-cell $v  \cdot  f \leq g  \cdot  u$. Then it is straightforward to check that 
\begin{equation}\label{oplaxstablediagonal}
    \xymatrix{
    D_v \armo[d]_{\sigma_v} \ar[dr]^{\sigma_g  \cdot  d} & \\
    B \prar[r]_{\delta} & C
    }
\end{equation}
is an oplax diagonal morphism for the oplax square through the 2-cells in \eqref{stablediag}. We remark moreover that, if $\cat{C} $ is an ordinary category, then the lower triangle determined by the oplax weak orthogonality is indeed commutative, i.e. $ v = g  \cdot  \delta $.

Moreover if $f \in ^{\oplorth}\overline{\mor{M}}$, then $f \oplorth m_f$ which implies that there exists an oplax diagonal morphism $\delta$ for the commutative square $ \ton{e_f , \id{B}} : f \longrightarrow m_f$. Such oplax diagonal morphism is total due to the properties mentioned in \ref{orderproperties}. Writing explicitly the diagrams it is easy to prove that $\varphi_{\delta}$ is an oplax diagonal morphism for $\ton{e_{\varphi_f}, \id{B}} : \varphi_f \longrightarrow m_{\varphi_f}$ in $\cat{C}$. One can prove by simple calculations that this yields that $\varphi_f \oplorth \mor{M}$. Henceforth $f$ belongs to $\overline{\mor{E}}$.

In a similar fashion we can prove that for any partial map $g \in ^{\oplorth} \overline{\mor{E}}$ there exists an oplax diagonal filler $\delta'$ for the square $\ton{\id{A} , m_g} : e_g \longrightarrow g$. Again writing explicitly the compositions one can prove that $\varphi_{\delta '} = \sigma_g  \cdot  \Tilde{\varphi_{\delta '}}$ for some $\Tilde{\varphi_{\delta '}}$ that is the oplax diagonal of the oplax square $\ton{\id{A} , m_{\varphi_g}} : e_{\varphi_g} \longrightarrow \varphi_g$. This yields that $\varphi_g$ lies in $\mor{M}$, hence $g \in \overline{\mor{M}}$.

In conclusion, $\ton{\overline{\mor{E}},\overline{\mor{M}}}$ is a lax weak factorisation system.

\subsection{Functoriality from total factorisation systems to oplax \textsc{wfs}}
Now we aim to prove that if the $\Ord$-enrichment on $\cat{C}$ is discrete, then the property of being functorial is carried from the factorisation system on $\cat{C}$ to the one on $\Pcat{\cat{C}}$.

\begin{proposition}
Let $\cat{C}$ be an ordinary category,  $\ton{F,L,R,K}$ a functorial factorisation and $\ton{\mor{E},\mor{M}}$ a stable \textsc{wfs} underlying it. Then $\ton{\overline{\mor{E}},\overline{\mor{M}}}$ underlies an oplax functorial factorisation $\ton{\overline{F}, \overline{L}, \overline{R}, \overline{K}}$.
\end{proposition}

\proof
We can rewrite the factorisation \eqref{EMfactdiag} substituting $L\varphi_f$ and $R\varphi_f$ by $m_f$ and $e_f$.
We obtain the assignments $\overline{L}f = \ton{\sigma_{f}, L\varphi_f} $, $\overline{R}f = \ton{\id{K\varphi_f}, R\varphi_f }$ and $\overline{K}f = K \varphi_f $. Then our goal is to prove their functoriality.
We consider two composable oplax squares
\begin{equation}\label{oplaxcomposition}
    \vcenter{
    \xymatrix{
    A \prar[r]^u \ar@{{}{ }{}}@/_2.5mm/[dr]^[@!45]{\preceq} \prar[d]_{f} & C \prar[d]^g\prar[r]^{u'} \ar@{{}{ }{}}@/_2.5mm/[dr]^[@!45]{\preceq} &E \prar[d]^h \\
    B \prar[r]_{v} & D \prar[r]_{v'} & G.
    }
    }
\end{equation}
We need to prove that $\overline{K} \ton{u' \cdot u , v' \cdot v} =\overline{K} \ton{u' , v'} \cdot \overline{K} \ton{ u , v} $.

First we reproduce the process depicted in \eqref{stablediag}. We choose among the possible diagonal liftings $k = K\ton{ \varphi_{\overline{L}g \cdot u} \cdot o, \varphi_{v \cdot \overline{R}f}}$. This diagonal lifting fills the following commutative diagram
\begin{equation}\label{kdiago}
    \vcenter{
    \xymatrix{
    D_{v \cdot f} \armo[r]^{o} \ar[d]_{\sigma^{*}_{v \cdot \overline{R}f} \ton{\varphi_{\overline{L}f}}} & D_{g \cdot u} \ar[r]^{\varphi_{\overline{L}g \cdot u}} & \overline{K}g \ar[d]^{R\varphi_{g}} \\
    D_{v \cdot \overline{R}f} \ar[urr]|{k} \ar[rr]_{\varphi_{v \cdot \overline{R}f}}&& D.
    }
    }
\end{equation}
Then $\overline{K}\ton{u,v}= \ton{R\varphi_f {}^* \ton{\sigma_v} , k}$. Similarly we have $\overline{K}\ton{u',v'}= \ton{R\varphi_g {}^* \ton{\sigma_{v'}} , k'}$, where $k' = K\ton{\sigma _{\overline{R}h} {}^{*} \ton{\varphi_{\overline{L}h \cdot u'} \cdot o'} , \varphi_{v' \cdot \overline{R}g}}$.
Finally we have that $\overline{K}\ton{u' \cdot u , v' \cdot v}= \ton{\sigma_{v' \cdot v \cdot \overline{R}f} , k''}$. We point out that
\begin{equation}
    \sigma_{v' \cdot v \cdot \overline{R}f} = R \varphi_f {}^{*} \cdot \qua{\ton{\varphi_v \cdot \sigma_v {}^{*} \ton{R\varphi_f}}^{*} \ton{\sigma_{v'}}}
\end{equation}
due to the pullback properties of subsequent composition of partial maps. Moreover $k''$ is the diagonal morphism  chosen through $K$ for the square

\begin{equation}\label{k2diago}
    \vcenter{
    \xymatrix{
    D_{v' \cdot v \cdot f} \armo[r]^{o * o'} \ar[d]_{\sigma^{*}_{v' \cdot v \cdot \overline{R}f} \ton{\varphi_{\overline{L}f}}} & D_{h \cdot u' \cdot u} \ar[r]^{\varphi_{\overline{L}h \cdot u' \cdot u }} & \overline{K}h \ar[d]^{R\varphi_{h}} \\
    D_{v' \cdot v \cdot \overline{R}f} \ar[urr]|{k''} \ar[rr]_{\varphi_{v' \cdot v \cdot \overline{R}f}}&& F.
    }
    }
\end{equation}

We write explicitly the composition $\overline{K} \ton{u' , v'} \cdot \overline{K} \ton{ u , v} $
\begin{equation}\label{Kbigcompo}
    \vcenter{
    \xymatrix{
    P  \armo[d]_{k^{*}\ton{R\varphi_{g}{}^{*}\ton{\sigma_{v'}}}}  \ar[dr]^{\ton{R\varphi_{g}{}^{*}\ton{\sigma_{v'}}}^*\ton{k}} && \\
    D_{v \cdot \overline{R}f} \armo[d]_{R\varphi_{f}{}^{*}\ton{\sigma_v}} \ar[dr]|{k} & D_{v' \cdot \overline{R}g} \armo[d]|{R\varphi_{g}{}^{*}\ton{\sigma_{v'}}} \ar[dr]^{k'} & \\
    \overline{K}f \prar[r]_{\overline{K} \ton{u , v}} & \overline{K}g \prar[r]_{\overline{K} \ton{u' , v'}} & \overline{K}h
    }
    }
\end{equation}

From \eqref{kdiago} we know that $\varphi_{v \cdot \overline{R}f} = \varphi_v \cdot \sigma_v{}^{*} \ton{R \varphi_f} = R\varphi_g \cdot k$.
Hence we have the equality of the two domains
\begin{equation}
    R\varphi_f{}^{*} \ton{\sigma_v} \cdot k^{*} \ton{R\varphi_g {}^{*} \ton{\sigma_{v'}}} = 
    R\varphi_f{}^{*} \cdot \qua{\ton{R \varphi_g \cdot k}^{*} \ton{\sigma_{v'}} } =
    R\varphi_f{}^{*} \cdot \qua{\ton{\varphi_v \cdot \sigma_v{}^{*}\ton{R\varphi_f}}^{*} \ton{\sigma_{v'}} }.
\end{equation}
Thus we can write $P= D_{v' \cdot v \cdot \overline{R}f}$.

Now we consider the second component of the partial maps.
We take the following diagram
\begin{equation}
    \vcenter{
    \xymatrix@C=4.5pc{
    D_{v' \cdot v \cdot f} \ar@{{}{ }{}}[drr]|{\fbox{1}} \armo[r]^{v' *o} \ar[d]|{\ton{ \sigma_{v' \cdot v \cdot \overline{R}f}}^* \ton{L\varphi_f}} & D_{v' \cdot g \cdot u} \ar[r]^{\sigma_{v' \cdot g} {}^{*} \ton{\varphi_u}} & D_{v' \cdot g} \ar@{{}{ }{}}[drr]|{\fbox{2}} \armo[r]^{o'} \ar[d]|{\ton{R\varphi_g {}^{*}  \ton{\sigma_{v'}}}^{*} \ton{L\varphi_g}}  & D_{h \cdot u'} \ar[r]^{\varphi_{\overline{L}h \cdot u' }} & \overline{K}h \ar[d]^{R\varphi_h}  \\
    D_{v' \cdot v \cdot \overline{R}f} \ar[rr]_{\ton{R\varphi_g {}^{*} \ton{\sigma_{v'}}}^* \ton{k}} & & D_{v' \cdot \overline{R}g} \ar[r]_{\sigma_{v'}{}^{*}\ton{R\varphi_g}} & D_{v'} \ar[r]_{\varphi_{v'}} & F.
    }
    }
\end{equation}

Due to the properties of pullbacks and the definition of each morphism we can prove that the diagram and its subsquares are actually commutative. In particular we have that \fbox{2} is the square used to deduce $k'$. We also point out that $\ton{ \sigma_{v' \cdot v \cdot \overline{R}f}}^* \ton{L\varphi_f}$ and $ \ton{R\varphi_g {}^{*} \ton{\sigma_{v'}}}^{*} \ton{L\varphi_g} $ are morphisms in $\mor{E}$, therefore factorising \fbox{1} we obtain 
\begin{equation}
    K\ton{ \sigma_{v' \cdot g} {}^{*} \ton{\varphi_u} \cdot \ton{v' *o}, \ton{R\varphi_g {}^{*} \ton{\sigma_{v'}}}^* \ton{k}} = \ton{R\varphi_g {}^{*} \ton{\sigma_{v'}}}^* \ton{k}.
\end{equation}

Examining the definition of the arrows, we can prove that the outer square is exactly the square that is used to define $k''$ as a diagonal filler. Thus, since $k''$ is chosen through $K$, by functoriality we can conclude that
\begin{equation}
    k''= k' \cdot \ton{R\varphi_g {}^{*} \ton{\sigma_{v'}}}^* \ton{k}.
\end{equation}
 This was the last information needed to conclude that $\ton{\overline{F},\overline{L},\overline{R},\overline{K}}$ is an oplax functorial factorisation system.
\endproof

\remark  It is straightforward to verify that for a functorial \textsc{wfs} $\ton{F,L,R,K}$, if $\ton{R, \eta}$ and $\ton{L, \varepsilon}$ are part of a monad and a comonad, then their oplax correspondents carry oplax monadic and comonadic structures. In fact, left and right components work in the same way as for total maps and satisfy the same axioms. The right component is always total and bears no difference from the total case. The left component operates similarly and one only has to take into account partial domains. This applies as well for the distributivity laws that define \textsc{awfs}.

In conclusion, an \textsc{awfs} $\ton{\nat{L} , \nat{R}}$ on $\cat{C}$ induces an oplax \textsc{awfs} $\ton{\overline{\nat{L}} , \overline{\nat{R}}}$ on $\Pcat{\cat{C}}$.

\example \label{ExampleSet} We consider $\Pcat{\Set}$. We know that in \Set the two classes \textup{\Epi} and \textup{\Mono} are stable under pullback. We know as well that $\ton{\textup{\Mono} , \textup{\Epi}}$ is a stable weak factorisation system. This yields that $\ton{ \overline{\textup{\Mono}} , \overline{\textup{\Epi}}}$ is an oplax weak factorisation system for $\Pcat{\Set}$. We also have that $\ton{\textup{\Epi} , \textup{\Mono}}$ is a stable orthogonal factorisation system. This yields that $\ton{ \overline{\textup{\Epi}} , \overline{\textup{\Mono}}}$ in an oplax \textsc{awfs}.

\begin{remark}
We conclude pointing out that $\ton{\overline{\Mono}, \overline{\Epi}}$-factorisations are not unique. In fact given a partial map $f\neq \emp$, it can be factorised as
\begin{equation}
    \vcenter{\xymatrix@C=5pc@R=3pc{
    D_f \armo[d]_{\sigma_f} \ar[dr]|{\id{D_f} \times \varphi_f} & D_f \times B \armo[d]|{\id{D_f \times B}} \ar[dr]^{\pi_B} & \\
    A \prar[r] & D_f \times B \prar[r] & B,
    }}
\end{equation}
but it can be also factorised as
\begin{equation}
    \vcenter{\xymatrix@C=5pc@R=3pc{
    D_f \armo[d]_{\sigma_f} \ar[dr]|{i_{D_f}} & D_f \amalg B \armo[d]|{\id{D_f \amalg B}} \ar[dr]^{f \amalg \id{B}} & \\
    A \prar[r] & D_f \amalg B \prar[r] & B.
    }}
\end{equation}
\end{remark}

\subsection{From Oplax \textsc{wfs} on partial maps to \textsc{wfs} on total maps}
We consider a category of partial maps $\Pcat{\cat{C}}$ and $\ton{\mor{L},\mor{R}}$ either a lax or an oplax weak factorisation system on $\Pcat{\cat{C}}$. We would like to analyse what kind of structure it generates on $\cat{C}$. 

We start by considering whether the orthogonality relations are preserved through this restriction. If we restrict to total maps, then lax and oplax squares reduce to commutative ones. We aim to prove that two total maps that are oplax weakly orthogonal ($f \oplorth g$), then they are also \textit{weakly orthogonal} ($f \boxslash g$) in the ordinary sense.

\begin{lemma}
Let $\ton{\mor{L},\mor{R}}$ be an oplax weak factorisation system. Then $\ton{\mor{L} \cap \textup{\textbf{Tot}}} \boxslash \ton{\mor{R} \cap \textup{\textbf{Tot}}}$.
\end{lemma}
\proof
We consider a commutative square formed by total maps and such that $l \in \mor{L}$ and $r \in \mor{R}$; then there exists a partial map $\delta$ that is an (op)lax diagonal morphism as
\begin{equation}\label{genlift}
    \vcenter{
    \xymatrix{
    A \ar[d]_l \ar[r]^u & C \ar[d]^r \\
    B \ar[r]_v \ar[ur]^{\delta} & D.
    }}
\end{equation}

If $\ton{\mor{L},\mor{R}}$ is an oplax weak factorisation system, then we have that $v \preceq r  \cdot  \delta$ and being $v$ a total map, it is an equality. Moreover $r  \cdot  \delta$ is then total and, by the rules of composition, $\delta$ is total as well. What said yields that the upper triangle must be commutative as well and therefore the oplax weak orthogonality relation restricts to a strict weak orthogonal relation among total maps.
\endproof

\begin{proposition}
Any total morphism admits an $\ton{\mor{L},\mor{R}}$-factorisation composed of total morphisms as well.
\end{proposition}
\proof
For every partial map $f$ and its $\ton{\mor{L},\mor{R}}$-factorisation we want to obtain another factorisation with a total right component, since we already know that, if $f$ is total, then its $\mor{L}$-component has to be total. We can build the following diagram for any partial morphism $f$:
\begin{align*}
    &\vcenter{\xymatrix{
    D_f \armo[d]_{\sigma_f} \ar[dr]^{\varphi_f} & \\
    A \prar[r]_f & B
    }}
    &
    &\longmapsto
    &
    &\vcenter{\xymatrix{
    D_f \armo[d] \ar[dr]^{\overline{\varphi_l}} \armo@/_2pc/[dd]_{\sigma_f} \ar@/^2.5pc/[ddrr]^{\varphi_f} && \\
    D_{l} \armo[d]_{\sigma_l} \ar[dr]^{\varphi_l} & D_{r} \armo[d]|{\sigma_r} \ar[dr]^{\varphi_r} & \\
    A \prar[r]_{l_f} & K \prar[r]_{r_f} & B
    }}
    &
    &\longmapsto
    &
    &\vcenter{\xymatrix{
    D_f \armo[d]|{\id{D_f}} \ar[dr]^{\overline{\varphi_l}} \armo@/_2pc/[dd]_{\sigma_f} \ar@/^2.5pc/[ddrr]^{\varphi_f} && \\
    D_f \armo[d]_{\sigma_f} \ar[dr]^{\overline{\varphi_l}} & D_{r} \armo[d]|{\id{D_r}} \ar[dr]^{\varphi_r} & \\
    A \prar[r]_{\overline{l_f}} & D_{r} \prar[r]_{\overline{r_f}} & B.
    }}
\end{align*}

Now we must understand whether $\overline{l_f} \in \mor{L}$ and $\overline{r_f} \in \mor{R}$. We first prove this helpful fact.

\begin{remark}
We consider the following adjoint morphisms $\ton{\id{D_{r}} , \sigma_r} \dashv \ton{ \sigma_r , \id{D_{r}}}$. We define $\overline{l_f} = \nu  \cdot  l_f $ and $\overline{r_f} = r_f  \cdot  \mu $ and by the adjunction 2-cells we obtain that
\begin{equation*}
    \begin{cases}
    \mu  \cdot  \overline{l_f} = \mu  \cdot  \nu  \cdot  l_f \preceq l_f \\
    \overline{r_f}  \cdot  \nu = r_f  \cdot  \mu  \cdot  \nu \preceq r_f.
    \end{cases}
\end{equation*}
By directly computing we actually get that $r_f = \overline{r_f}  \cdot  \nu$. Moreover $D_l =D_f$ if and only if $l_f = \mu  \cdot  \overline{l_f}$.
\end{remark}

We consider $\overline{r_f}$ and we prove that it lies in $\mor{R}$. We take $l \in \mor{L}$ and an oplax square $\ton{u , v}: l \longrightarrow \overline{r_f}$. We can build the diagram
\begin{equation*}
    \xymatrix@R=1pc@C=1pc{
    C \ar[ddd]_l \ar[rrr]^{u} \ar[drr]^{\mu  \cdot  \nu} &&& D_r \ar[ddd]^{\overline{r_f}} \ar[ld]|{\mu} \\
    && K \ar[ddr]_{r_f} & \\
    &&& \\
    D \ar[rrr]_{v} \ar@{-->}[rruu]|{\delta} &&& B
    }
\end{equation*}
where the upper and right triangles are commutative and $\delta$ is an oplax diagonal morphism lifting $l \in \mor{L}$ against $r_f \in \mor{R}$.
We consider a diagonal $\nu  \cdot  \delta$. Then considering the outer square we have
\begin{equation*}
    \begin{cases}
    \nu  \cdot  \delta  \cdot  l \succeq \nu  \cdot  \mu  \cdot  u = u ;\\
    \overline{r_f}  \cdot  \nu  \cdot  \delta = r_f  \cdot  \delta \preceq v.
    \end{cases}
\end{equation*}
Hence $\nu  \cdot  \delta$ is an oplax diagonal morphism for $l$ against $\overline{r_f}$ and in conclusion $\overline{r_f} \in \mor{R}$.

Let us consider $\overline{l_f}$. For any $r \in \mor{R}$ part of an oplax square as the outer diagram

\begin{equation*}
    \xymatrix@R=1pc@C=1pc{
    A \ar[ddd]_{\overline{l_f}} \ar[rrr]^{u} \ar[ddr]^{l_f} &&& C \ar[ddd]^{r}  \\
    &&& \\
    & K \ar@{-->}[rruu]|{\delta} \ar[drr]^{\nu  \cdot  \mu} \ar[dl]|{\nu} && \\
    D_r \ar[rrr]_{v}  &&& D,
    }
\end{equation*}
again the lower and left triangles are commutative and $\delta$ is a lax/oplax diagonal morphism lifting $l_f \in \mor{L}$ against $r \in \mor{R}$. Now the diagonal morphism for the outer diagram is $\delta  \cdot  \mu$ and the proof proceeds analogously. Hence we have that $\overline{l_f} \in \mor{L}$.
\endproof

\begin{remark}
As for the lax case regarding \eqref{genlift}, we obtain $u \preceq \delta  \cdot  l$ and $r  \cdot  \delta \preceq v$, which do not imply in general that $\delta $ is total. Hence we cannot deduce the commutativity of any triangle in the diagram. Nonetheless one can reproduce the same process of extracting a total factorisation from any factorisation for total maps. This process in the lax case is in fact successful, but only for morphisms $f$ such that $D_l=D_f$.
\end{remark}

We conclude this section by proving the following proposition.
\begin{proposition}
Let $\Pcat{\cat{C}}$ be a category of partial maps and $\ton{\mor{L}, \mor{R}}$ an oplax \textsc{wfs}. Then
$\ton{\mor{L}_{\textup{\textbf{Tot}}} , \mor{R}_{\textup{\textbf{Tot}}}}$ is a weak factorisation system for $\cat{C}$.
\end{proposition}
\proof
As we have seen above the oplax weak orthogonality relation restricts to a weak orthogonality relation among total maps, henceforth $\mor{L}_{\textbf{Tot}} \boxslash \mor{R}_{\textbf{Tot}}$.
If $f \boxslash \mor{R}_{\textbf{Tot}}$, then $f$ has a $\ton{\mor{L}_{\textbf{Tot}} , \mor{R}_{\textbf{Tot}}}$-factorisation $f = r_f  \cdot  l_f$. The commutative square $\ton{l_f , \id{B}}: f \longrightarrow r_f$ has a diagonal morphism $\delta$ and, in particular, $r_f  \cdot  \delta = \id{B}$. Since $r_f$ is a coretract and $\mor{L}$ is closed under coretract composition to the left, then it is easily proved that $f \in \mor{L}_{\textbf{Tot}}$, thus $^{\boxslash}\mor{R}_{\textbf{Tot}} = \mor{L}_{\textbf{Tot}}$. A dual argument proves, using retract closure of $\mor{R}$, that $\mor{L}_{\textbf{Tot}}^{\boxslash} = \mor{R}_{\textbf{Tot}}$. This yields the thesis.
\endproof

\section{Factorisations for pointed categories of partial maps}
Our goal for the following section is to discuss a process to obtain lax and oplax weak factorisation systems for pointed categories of partial maps. In the first part we will show how this process is carried out for $\Pcat{\Set}$, and then we will try to generalise it for any pointed category of partial maps.

\subsection{\Set with partial maps}
In the following section we will consider the category of partial maps among sets, together with the $\Ord$-enrichment induced by the discrete order in \Set.
In $\Pcat{\Set}$ $\emp$ is a zero object. In fact for every pair of sets $A,B$ we have the zero map
\begin{equation}
    \vcenter{
    \xymatrix{
    \emp \armo[d] \armo[dr] & \\
    A \prar[r]_{\emp_{A,B}} & B.
    }
    }
\end{equation}
We define the class of morphisms
\begin{equation}
\mor{O} = \bra{\emp_{A,B} \vert A,B \in \Set}.
\end{equation}
We will apply Proposition \ref{smallobjectlax} and its dual to $\mor{O}$.
\begin{itemize}
    \item First we consider an $f \in {}^{\lorth}\mor{O}$. There exists a lax diagonal morphism $f^{*}$ for the lax square $\ton{\id{A} , \id{B}} : f \longrightarrow \emp_{A,B}$ and hence there exists the 2-cell $\id{A} \preceq f^{*} \cdot f $. We notice that the existence of this 2-cell yields that $f$ is a left adjoint to $f^{*}$, hence the lax factorisation system is trivial.

    \item We consider now $f \in \mor{O}^{\lorth}$. We notice that any lax square $\ton{u,v} : \emp \longrightarrow f $ yields that $f \cdot u = \emp$. It is straightforward to prove that there exists a lax diagonal lifting if and only if $u = \emp$. We deduce that
    \begin{equation}
    \mor{O}^{\lorth} = \bra{f \vert f \cdot u = \emp \: \Rightarrow \: u = \emp}.
    \end{equation}
    This in $\Pcat{\Set}$ is equivalent to have that $f$ is a total map. In fact, any total map belongs to $\mor{O}^{\lorth}$ trivially. Moreover if $f: A \rightarrow B$ is not total, then any non-zero map $g$ whose image is a subset of $A \backslash D_{f}$, is $\emp $ when composed with $f$. We conclude that $\mor{O}^{\lorth} = \textbf{\textup{Tot}}$. As seen in Section \ref{TotDom}, the factorisation system that is generated is therefore $\ton{\overline{\mor{S}}, \textbf{Tot}}$, since in $\Pcat{\Set}$ we have that $\textbf{PLA} = \overline{\mor{S}}$.

    \item Then we consider ${}^{\oplorth}\mor{O}$. Similarly to the previous case, we can deduce that $f \in {}^{\oplorth}\mor{O}$ if and only if for any morphism $v$, $v \cdot f = \emp$ implies that $v= \emp$, and it is equivalent to surjectivity of $f$. Hence ${}^{\oplorth}\mor{O} = \overline{\Epi}$. As discussed before in Example \ref{ExampleSet} $\overline{\Epi }$ is part of the oplax \textsc{awfs} $\ton{\overline{\Epi} , \overline{\Mono}}$.

    \item Finally we consider $f \in \mor{O}^{\oplorth}$. Then in particular there exists an oplax diagonal lifting $f_{*}$ for the oplax square $\ton{\id{A}, \id{B}} : \emp \longrightarrow f$ and the 2-cell $\id{A} \preceq f \cdot f_{*}$. Moreover the existence of such morphism and 2-cell is easily shown to be a sufficient condition that implies $f \in \mor{O}^{\oplorth}$, since $f_{*}$ is a tool to build oplax diagonal morphisms for any other oplax squares. We observe that, in the context of $\Pcat{\Set}$, this condition is equivalent to $f \in \overline{\Epi}$. Again we know that $\overline{\Epi}$ is part of the oplax factorisation system $\ton{\overline{\Mono} , \overline{\Epi}}$, as shown in Example \ref{ExampleSet}.
\end{itemize}

\subsection{Factorisations for pointed categories of partial maps}
In the following subsection, we will try to expand the process described for $\Pcat{\Set}$ and $\mor{O}$ and apply it to other categories with similar properties.

Along this section we will assume that \textit{the $\Ord$-category $\cat{C}$ has an initial object $I$ and that all initial morphisms $i_X$ lie in $\mor{S}$}, implying in particular that they are monomorphisms.
We notice that $I$ is still an initial object in $\Pcat{\cat{C}}$ and that for any $A,B$ the partial morphism $\iota = \ton{i_A,i_B}$ is minimal in $\cat{C}\ton{A,B}$. In fact for any partial map $f$, the arrow $i_{D_f}$ shows that $\iota \leq f$. On the other hand if $f \leq \iota$, then there exists $s: D_f \rightarrowtail I$, where $s$ is monomorphic, which yields that $f=\iota$.
We can consider now the class of minimal maps
\begin{equation}
    \mor{O}= \bra{\iota= \ton{i_A,i_B}: A \longrightarrow B \vert A;B \in \textbf{Ob}\ton{\cat{C}}}.
\end{equation}

\begin{remark}\label{factzeroobj}
We observe that in $\Pcat{\cat{C}}$ the initial object $I$ of $\cat{C}$ is a zero object whenever $I$ is either a zero object or a strict initial object. 
\end{remark}
This is true since the choice for the component $\varphi_-$ becomes unique when the codomain is $I$ under the said assumption. We recall that cartesian closed categories, such as $\Set$, \textbf{Cat}, any topos, and  distributive categories have strict initial objects.

We remark that the property described in Remark \ref{factzeroobj} is not always needed. In fact, we are interested in the property of minimal maps of being \textit{left} or \textit{right absorbent}, meaning that when a minimal map is composed on the left or on the right with any other composable morphism, then the composition is a minimal morphism in the corresponding Hom-Set.

\begin{lemma}\label{factabsorbent}
For any $\Pcat{\cat{C}}$, a minimal map $f$ such that $D_f=I$ is right absorbent. Whenever $I$ is actually a zero-object, then it is both left and right absorbent.
\end{lemma}
This is trivial considering that the partial domain of the composition is a subobject of the partial domain of the first morphism and $I$ admits only itself as a subobject.

Henceforth the hypothesis that $I$ is a zero-object in $\Pcat{\cat{C}}$ is needed only while discussing  the left complements $^{\lorth}\mor{O}$ and $^{\oplorth}\mor{O}$.

We will now apply Proposition \ref{smallobjectlax} to the class $\mor{O}$; this result will enable us to reproduce the process in $\mor{O}$ in two directions for the lax case.

\begin{enumerate}
    \item We consider $^{\lorth} \mor{O}$. Let  $f: A \longrightarrow B$ be a partial map laxly weakly orthogonal to $\mor{O}$. Then there is a lax diagonal $f^*$ for the lax square $\ton{\id{A},\id{B}} : f \longrightarrow 0$. While $0 \cdot f^* = 0 \leq \id{B}$ always exists by minimality of zero maps, the 2-cell $\id{A} \leq f^* \cdot f$ does not exist in general. This property of $f$ of having a paired arrow $f^{*}$ such that $\id{A} \leq f^* \cdot f$ is also a sufficient condition for $f$ to be in $^{\lorth} \mor{O}$.
    In conclusion
    \begin{equation}\label{Uclass}
        \mor{U} = ^{\lorth} \mor{O}= \bra{f \textup{ } \vert \textup{ } \id{A} \leq f^* \cdot f \textup{ for some } f}.
    \end{equation}
    Unfortunately we could not find a general description of $\mor{U}^{\lorth}$. Still we present our:
    \begin{conjecture}
    The complement $\mor{U}^{\lorth } = \bra{f\vert  f \cdot f^* \leq \id{B} \textup{ for some } f^*}$, the intersection $\mor{U}^{\lorth} \cap \mor{U}$ being exactly the left adjoint morphisms.
    \end{conjecture}
\end{enumerate}

\begin{remark}
If the partial order on $\Pcat{\cat{C}}$ is induced by the discrete partial order on $\cat{C}$, then we can refine the description of some complements of $\mor{O}$.
In fact, this $\Ord$-enrichment yields that $^{\lorth}\mor{O}$ is actually the class of coretract partial morphisms $\textbf{LI}$, which is the class of left adjoint morphisms. In this case the factorisation produced is the trivial $\ton{\textbf{LA}, \textbf{All}}$.
\end{remark}

\begin{enumerate}
    \item[2.] We now consider $\mor{O}^{\lorth} $. If $\ton{u,v}: 0 \longrightarrow f$ is a lax square, then $f \cdot u=0$, since zero maps are absorbent and minimal. If $f \in \mor{O}^{\lorth}$, then there exists a lax diagonal morphism $\delta$ as in
    \begin{equation}\vcenter{
    \xymatrix{
    A \ar[r]^{u} \ar@{{}{ }{}}@/^3mm/[d]^(.35)[left]{\geq} \ar[d]_{0} & C \ar[d]^{f} \\
    B \ar[ur]|{\delta} \ar[r]_{v} & D \ar@{{}{ }{}}@/^3mm/[u]^(.35)[left]{\geq}.
    }}
    \end{equation}
    We notice that a necessary and sufficient condition for the existence of such $\delta$ is that $u=0$. So we have that 
    \begin{equation}\label{densedomain}
    \mor{O}^{\lorth} =  \mor{DD} =\bra{f\vert f \cdot u = \iota \Longrightarrow u=\iota }.
    \end{equation}
    The notation $\mor{DD}$ denotes those maps whose partial domain is maximal as a proper S-subobject. Inspired by the example of topological spaces below, we chose to call \textit{dense domain partial maps}.
\end{enumerate}

\begin{remark}
As we proved above, in \Set we have that $\mor{DD}=\textup{\textbf{Tot}}$, and therefore $\ton{^{\lorth}\mor{DD}, \mor{DD} } = \ton{\overline{\mor{S}}, \textup{\textbf{Tot}}}$. In general we know that $\mor{DD} \supseteq \textup{\textbf{Tot}}$, but the other inclusion is not always true. In fact, we have counterexamples of partial morphisms which are not total, but have a dense domain:
\begin{itemize}
    \item in $\san{Ab}$ maps such as
    \begin{align*}
        &\xymatrix{
        \nat{Z} \ar@{>->}[d]_{2} \ar[dr] && \\
        \nat{Z} \prar[r] & \nat{Z}
        }
        &
        &\xymatrix{
        \nat{Z} \ar@{>->}[d]_{i} \ar[dr] && \\
        \nat{Q} \prar[r] & \nat{Z}
        }
    \end{align*}
    are not total, but it is easily proved that they have dense domains.
    \item \label{densemaps} for the category of topological spaces equipped open maps, we have that a domain is dense exactly when the domain is a topologically dense subobject of the domain, so any morphism $f$ such that $\sigma_f=j: \left[ 0,1 \right[ \longrightarrow \qua{0,1}$ is not total and yet it has a dense domain.
\end{itemize}

Since total maps have always a dense domain, we have that $\ton{\overline{\mor{S}},\textup{\textbf{Tot}}} \leq \ton{^{\lorth}\mor{DD}, \mor{DD}}$.
\end{remark}

We briefly state the two counterparts for $\Coplax$ that arise in a similar fashion.

\begin{enumerate}
    \item The left oplax complement is
    \begin{equation}\label{denseimage}
        \mor{DI} = ^{\oplorth} \mor{O} = \bra{ f\vert v \cdot f =\iota \Longrightarrow v=\iota },
    \end{equation}
    where the notation $\mor{DI}$ refers to dense image maps among partial maps.
    
    \item Considering $^{\oplorth} \mor{O}$ one can prove that
    \begin{equation}\label{Vclass}
        \mor{V} = \mor{O}^{\oplorth} = \bra{f\vert \id{B} \leq f \cdot f_* \textup{ for some } f_*}.
    \end{equation}
    We conjecture again:
    \begin{conjecture}
    The complement $\mor{V}^{\oplorth } = \bra{f\vert   f \cdot f_* \leq \id{A} \textup{ for some } f_*}$, the intersection being exactly the right adjoint morphisms.
    \end{conjecture}
\end{enumerate}

\begin{remark}
We consider now an $f \in \mor{O}^{\oplorth}$; then there exists an oplax diagonal morphism $\delta$ for the oplax square $\ton{\id{A} , \id{B}}: \iota_{A,B} \longrightarrow f$. Since the identity is total, the lower triangle is indeed commutative, therefore $f$ is a split epimorphism (\textsc{li}). On the other hand if $f$ has a right inverse $f'$, then in any oplax square determined by the cell $v  \cdot  \iota \leq f  \cdot  v$, the morphism $f'  \cdot  v $ is an oplax diagonal morphism. This yields that $\mor{O}^{\oplorth} = \textbf{LI}$.

Considering the fact that the intersection of the two classes must be the class of right adjoint morphisms, i.e. $f$ such that $\varphi_f$ is an isomorphism, we conjecture the following.

\begin{conjecture}
The oplax weak orthogonal complement $^{\oplorth}\textup{\textbf{LI}}$ is $\overline{\mor{S}}$.
\end{conjecture}
\end{remark}

\begin{remark}
Looking carefully at the arguments we notice that the process described does not use specific tools for categories of partial maps. In fact the main ingredients is to consider an $\Ord$-enriched category such that every Hom-Set has a minimal element and the class of such minimal elements are absorbent. The description of the complements is essentially the same as above and can be carried out for pointed $\Ord$-enriched categories such that 0-maps are minimal in their Hom-Sets.
\end{remark}

Even if it has been difficult to give a better description for such complements, we remark that in general these classes appear to be non-trivial. To reduce to cases where these classes are trivial we actually have to impose strong restrictions on $\cat{C}$ and $\mor{S}$, such as having only split epimorphisms among the arrows or similar assumptions.

\section*{Conclusion}
In conclusion we have presented how the introduction of a lax weak orthogonality relation induces new notions of factorisation systems that carry similar facets to their discrete counterparts. We have defined the general notion of \textsc{lwfs} and \textsc{lffs}, showing how they relate to each other and discussed examples and constructions in categories of partial maps.

Future developments of this study will be expanding the set of examples and applications of these structures.

Moreover an interesting future direction is to explore the connection between this work and the characterization of lax orthogonal factorisation systems presented by John Bourke and Charles Walker in \cite{walker2020characterization}, in particular studying the relations between \textit{down factorization systems}, that they introduce, and the lax and oplax factorisation systems that we have introduced; in particular some of the examples we provide that carry similar properties to their structures.

Our work on partial morphisms is also connected to the recent \cite{Cockett2020Latent}, in which some similar constructions are introduced. In particular this work focuses as well on the relation between factorisation systems on a category of partial maps and the stable factorisation systems on the base category.

\section*{Acknowledgment}
This work was done during the preparation of the author’s PhD thesis, under the supervision of Maria Manuel Clementino, whom the author thanks for proposing the investigation and advising the whole study.

\bibliographystyle{apalike}

\bibliography{mylibrary}

\end{document}